\documentclass[11pt,a4paper]{article}

\usepackage[T1]{fontenc}
\usepackage{epsfig}
\usepackage{amsmath, amsthm}
\usepackage{amsfonts,amssymb}

\newtheorem{conj}{Conjecture}[section]

\newtheorem{prop}[conj]{Proposition}
\newtheorem{coro}[conj]{Corollary}

\newcommand{\de}{\mathrm{d}}
\newcommand{\e}{\mathrm{e}}

\newcommand{\proj}{\mathrm{P}}

\newcommand{\R}{\mathbb{R}}

\newcommand{\eps}{\varepsilon}

\begin{document}

\title{On the improvement of concavity of convex measures}

\author{Arnaud Marsiglietti}

\date{}

\maketitle
{\raggedright Universit\'e Paris-Est, LAMA (UMR 8050), UPEMLV, UPEC, CNRS, F-77454, Marne-la-Vall\'ee, France} \\

\noindent
arnaud.marsiglietti@u-pem.fr \\

{\raggedright {\scriptsize The author was supported in part by the Agence Nationale de la Recherche, project GeMeCoD (ANR 2011 BS01 007 01).}}

\begin{abstract}
We prove that a general class of measures, which includes $\log$-concave measures, is $\frac{1}{n}$-concave according to the terminology of Borell, with additional assumptions on the measures or on the sets, such as symmetries. This generalizes results of Gardner and Zvavitch~\cite{GZ}.
\end{abstract}


\section{Introduction}

The classical Brunn-Minkowski inequality asserts that for all non-empty Borel sets $A,B \subset \R^n$ and for every $\lambda \in [0,1]$, one has
\begin{eqnarray}\label{BM}
|(1-\lambda) A + \lambda B|^{\frac{1}{n}} \geq (1-\lambda)|A|^{\frac{1}{n}} + \lambda |B|^{\frac{1}{n}},
\end{eqnarray}
where
$$ A + B = \{a + b; a \in A, b \in B\} $$
denotes the Minkowski sum of $A$ and $B$ and where $|\cdot|$ denotes the Lebesgue measure. 

The Brunn-Minkowski inequality is a beautiful and powerful inequality in Geometry and Analysis, leading to many interesting consequences. For more information on this inequality and its influences on several mathematical theories, see the survey by Gardner~\cite{Gardner}. See also the book by Schneider \cite{Schneider}, as a general reference in Convex Geometry.

Recently, Gardner and Zvavitch~\cite{GZ} proved that the Gaussian measure $\gamma_n$ on $\R^n$, defined by
$$ \de \gamma_n(x) = \frac{1}{(2 \pi)^{\frac{n}{2}}} \e^{-\frac{|x|^2}{2}} \, \de x, \quad x \in \R^n $$
where $|\cdot|$ denotes the Euclidean norm, satisfies a Brunn-Minkowski-type inequality of the form~(\ref{BM}) for some special classes of sets. More precisely, Gardner and Zvavitch showed that for coordinate boxes ($i.e.$ a product of intervals) $A, B \subset \R^n$ that contain the origin, or for $A,B \subset \R^n$ that are dilates of the same symmetric convex set, and for every $\lambda \in [0,1]$, one has
\begin{eqnarray}\label{GBM}
\gamma_n((1-\lambda) A + \lambda B)^{\frac{1}{n}} \geq (1-\lambda) \gamma_n(A)^{\frac{1}{n}} + \lambda \gamma_n(B)^{\frac{1}{n}}.
\end{eqnarray}
They also conjectured that inequality~(\ref{GBM}) holds for all convex symmetric sets $A,B \subset \R^n$.

As a consequence of the Pr\'ekopa-Leindler inequality \cite{Prekopa71}, \cite{Leindler}, \cite{Prekopa}, the Gaussian measure satisfies for all Borel sets $A,B \subset \R^n$ and for every $\lambda \in [0,1]$,
\begin{eqnarray}\label{G}
\gamma_n((1-\lambda) A + \lambda B) \geq \gamma_n(A)^{1-\lambda} \gamma_n(B)^{\lambda}.
\end{eqnarray}
Using the terminology of Borell~\cite{Borell} (see Section~2 below for further details), inequality~(\ref{G}) means that the Gaussian measure is a $\log$-concave measure. By comparing means, one sees that inequality~(\ref{GBM}) is stronger than inequality~(\ref{G}), hence the results of Gardner and Zvavitch improve the concavity of the Gaussian measure by showing that this measure is $\frac{1}{n}$-concave if restricted to special classes of sets. \\

We will see in this paper that the results of Gardner and Zvavitch can be extended to a more general class of measures called {\it convex measures} that includes the Gaussian measure. This is the mathematical underlying idea of the Gaussian Brunn-Minkowski inequality~(\ref{GBM}), $i.e.$ under symmetry assumptions one can improve a certain property; in this case it is the concavity of a measure. However, we will see that symmetries are not the only hypothesis that allow improvement of concavity properties of a measure. 

This paper is devoted to the study of the following question: \\

\paragraph*{{\bf Question 1.}}

For which value $s \in [-\infty, + \infty]$, for which class $\mathcal{M}$ of measures on $\R^n$ and for which class $\mathcal{C}$ of couples of Borel subsets of $\R^n$ one has
\begin{eqnarray}\label{question}
\mu((1-\lambda)A + \lambda B) \geq \left((1-\lambda) \mu(A)^s + \lambda \mu(B)^s\right)^{\frac{1}{s}}
\end{eqnarray}
for every $\mu \in \mathcal{M}$, for every $(A,B) \in \mathcal{C}$ such that $\mu(A)\mu(B) > 0$ and for every $\lambda \in [0,1]$?

The right-hand side of inequality~(\ref{question}) has to be interpreted by $\mu(A)^{1-\lambda} \mu(B)^{\lambda}$ for $s=0$, by $\min(\mu(A),\mu(B))$ for $s = -\infty$ and by $\max(\mu(A),\mu(B))$ for $s = +\infty$. \\

Borell~\cite{Borell} (see also~\cite{BL}) answered Question~1 when $\mathcal{M}$ is the class of $s$-concave measures on $\R^n$, $s \in [-\infty, + \infty]$, and when $\mathcal{C}$ is the class of all couples of Borel subsets of $\R^n$ (see Section~2).

Question~1 has been explored for $s = 1$, when restricted to the Lebesgue measure, by Bonnesen~\cite{BF} and is still being studied (see $e.g.$~\cite{HN}). \\

The main results of this paper are contained in the following theorem: \\

\paragraph*{{\bf Theorem 1.}}

\begin{enumerate}
\item Let $\mu$ be an unconditional $\log$-concave measure on $\R^n$ and let $A$ be an unconditional convex subset of $\R^n$. Then, for every $A_1, A_2 \in \{\alpha A; \alpha~>~0\}$ and for every $\lambda \in [0,1]$, we obtain
$$ \mu((1-\lambda)A_1 + \lambda A_2)^{\frac{1}{n}} \geq (1-\lambda) \mu(A_1)^{\frac{1}{n}} + \lambda \mu(A_2)^{\frac{1}{n}}. $$
\item Let $\mu_i$, $1 \leq i \leq n$, be measures with densities $\phi_i : \R \to \R_+$ such that $\phi_i$ are non-decreasing on $(-\infty; 0]$ and non-increasing on $[0; + \infty)$. Let $\mu$ be the product measure of $\mu_1, \dots, \mu_n$ and let $A,B \subset \R^n$ be the product of $n$ Borel subsets of $\R$ such that $0 \in A \cap B$. Then, for every $\lambda \in [0,1]$, we obtain
$$ \mu((1-\lambda) A + \lambda B)^{\frac{1}{n}} \geq (1-\lambda) \mu(A)^{\frac{1}{n}} + \lambda \mu(B)^{\frac{1}{n}}. $$
\end{enumerate}

In the next section, we introduce new terminology and in the third section, we prove Theorem~1. In the last section, we discuss how these results improve concavity properties of the (extended) parallel volume.

\section{Preliminaries}

We work in the Euclidean space $\mathbb{R}^n$, $n \geq 1$, equipped with the $\ell_2^n$ norm $|\cdot |$. The closed unit ball is denoted by $B_2^n$, the unit sphere by $\mathcal{S}^{n-1}$ and the canonical basis by $\{e_1, \cdots, e_n\}$. We also denote by $|\cdot |$ the Lebesgue measure on $\mathbb{R}^n$. For $u \in \mathcal{S}^{n-1}$, we denote by $u^{\perp}$ the hyperplane orthogonal to $u$. For non-empty sets $A,B$ in $\mathbb{R}^n$ we define their \textit{Minkowski sum}
$$ A+B = \{a + b ; \, a \in A, b \in B \}. $$
The {\it parallel volume} of a non-empty Borel set $A \subset \R^n$ is the function defined on $\R_+$ by $t \mapsto |A + tB_2^n|$.

A subset $A \subset \R^n$ is said to be {\it symmetric} if $A = -A$. A function $f : \R^n \to \R$ is said to be {\it unconditional} if for every $(x_1, \cdots, x_n) \in \R^n$ and for every $(\eps_1, \cdots,\eps_n) \in \{-1,1\}^n$, one has $f(\eps_1 x_1, \cdots, \eps_n x_n) = f(x_1, \cdots, x_n)$. A subset $A \subset \R^n$ is said to be {\it unconditional} if the indicator function of $A$, denoted by $1_A$, is unconditional. A measure with density function is said to be {\it symmetric} (resp. {\it unconditional}) if its density function is even (resp. unconditional).

Let us recall the terminology and the results of $s$-concave measures introduced by Borell in \cite{Borell}. One says that a measure $\mu$ on $\mathbb{R}^n$ is \textit{$s$-concave}, $s \in [- \infty, + \infty]$, if the inequality
\begin{eqnarray}\label{s-concave}
\mu((1 - \lambda)A + \lambda B) \geq ((1-\lambda)\mu(A)^s + \lambda \mu(B)^s)^{\frac{1}{s}}
\end{eqnarray}
holds for all Borel subsets $A,B \subset \mathbb{R}^n$ such that $\mu(A)\mu(B) > 0$ and for every $\lambda \in [0,1]$. The limit cases are interpreted by continuity, as mentioned in the introduction. The $0$-concave measures are also called {\it log-concave measures}.

Notice that an $s$-concave measure is $r$-concave for every $r \leq s$. Thus, every $s$-concave measure is $-\infty$-concave. The $-\infty$-concave measures are also called {\it convex measures}.

From inequality~(\ref{G}), the Gaussian measure is a $\log$-concave measure and as a consequence of the Brunn-Minkowski inequality~(\ref{BM}), the Lebesgue measure is a $\frac{1}{n}$-concave measure.

For every $s \in [-\infty, +\infty]$, Borell gave a complete description of $s$-concave measures. In particular, for $s \leq \frac{1}{n}$, Borell showed that every measure $\mu$ that is absolutely continuous with respect to the $n$-dimensional Lebesgue measure, is $s$-concave if and only if its density function is a $\gamma$-concave function, with
$$ \gamma = \frac{s}{1-sn}  \in [-\frac{1}{n}, + \infty], $$
where a function $f : \R^n \to \R_+$ is said to be \textit{$\gamma$-concave}, with $\gamma \in [-\infty, +\infty]$, if the inequality
$$ f((1 - \lambda)x + \lambda y) \geq ((1-\lambda)f(x)^{\gamma} + \lambda f(y)^{\gamma})^{\frac{1}{\gamma}} $$
holds for every $x,y \in \R^n$ such that $f(x)f(y) > 0$ and for every $\lambda \in [0,1]$. As for the $s$-concave measures, the limit cases are interpreted by continuity. Notice that a 1-concave function is concave on its support, that $f$ is a $- \infty$-concave function if and only if $f$ has convex level sets, and that $f$ is a $+ \infty$-concave function if and only if $f$ is constant on its support. \\

Measures with $-\infty$-concave density function are a natural generalization of convex measures. The results of Borell show that a measure that has a $\gamma$-concave density function with $\gamma < -\frac{1}{n}$, does not satisfy a concavity property of the form~(\ref{s-concave}) (but satisfies other forms of concavity~\cite{DU}). However, we will show that if restricted to a special class of sets, such measures are $\frac{1}{n}$-concave.

In this paper, we call {\it sub-convex measure} a measure that has a $-\infty$-concave density function. Notice that convex measures are sub-convex.

\section{Brunn-Minkowski type inequality for sub-convex measures}

In this section, we partially answer Question~1 by investigating possible improvements of the concavity of sub-convex measures. Gardner and Zvavitch~\cite{GZ} noticed in the case of the Gaussian measure, that the position of the sets $A$ and $B$ plays an important role. Indeed, since for $s$-concave probability measures $\mu$, with $s \leq 0$, the density function tends to $0$ at infinity and the support can be equal to $\R^n$, one can find sets $A$ and $B$ such that $A$ contains the origin and $\frac{A+B}{2}$ is far from the origin. Thus for $r>0$, the inequality
$$ \mu\left(\frac{A+B}{2}\right)^r \geq \frac{\mu(A)^r + \mu(B)^r}{2} $$
will not be satisfied. Hence, the position of sets $A$ and $B$ is an inherent constraint of the problem. Notice also that in the definition of $s$-concave measures, the condition $\mu(A) \mu(B) > 0$ is already a constraint on the position of $A$ and $B$ with respect to the support of $\mu$. \\

Notice that Question~1 has an answer for $s = + \infty$ if $\mathcal{M}$ is the class of convex measures and if $\mathcal{C}$ is the class of couples of Borel sets with the same measure. Indeed, one then has for every $\lambda \in [0,1]$,
$$ \mu((1-\lambda)A + \lambda B) \geq \inf(\mu(A), \mu(B)), $$
by definition. Since $\mu(A) = \mu(B)$, it follows that
$$ \mu((1-\lambda)A + \lambda B) \geq \mu(A) = \max(\mu(A), \mu(A)) = \max(\mu(A), \mu(B)). $$

Notice also that for every measure $\mu$ and for all Borel sets $A,B$ such that $A \subset B$, one has for every $\lambda \in [0,1]$,
$$ \mu((1-\lambda)A + \lambda B) \geq \min(\mu(A),\mu(B)), $$
since in this case one has, $(1-\lambda)A + \lambda B \supset (1-\lambda)A + \lambda A \supset A$.

\subsection{The case of symmetric measures and symmetric sets}

Under symmetry assumptions, the best concavity one can obtain is $\frac{1}{n}$ by considering, for example, the Lebesgue measure, which fulfills a lot of symmetries (unconditional), and two dilates of $B_2^n$ (which are unconditional). This was noticed by Gardner and Zvavitch~\cite{GZ} also for the Gaussian measure. \\

A sufficient condition to ensure that a measure $\mu$ on $\R^n$ is $\frac{1}{n}$-concave in the class of dilates of a fixed Borel set $A \subset \R^n$ is that the function $t \mapsto \mu(tA)$ is $\frac{1}{n}$-concave. The following proposition gives a sufficient condition for this.

\begin{prop}\label{equiv}

Let $\phi : \R^n \to \R_+$ be a measurable function such that for every $x \in \R^n$, the function $t \mapsto \phi(tx)$ is  non-increasing on $\R_+$. Let $\mu$ be a measure with density function $\phi$ and let $A$ be a Borel subset of $\R^n$ containing $0$. If the function $t \mapsto \mu(\e^t A)$ is $\log$-concave on $\R$, then the function $t \mapsto \mu(t A)$ is $\frac{1}{n}$-concave on $\R_+$.

\end{prop}

\begin{proof}
Let $\mu$ be a measure with density function $\phi$ satisfying the assumptions of Proposition~\ref{equiv} and let $A$ be a Borel subset of $\R^n$ containing $0$. Let us denote $F(t) = \mu(tA)$, for $t \in \R_+$. Notice that $F$ is non-decreasing and continuous on $\R_+$. By assumption, the function $t \mapsto F(\e^t)$ is $\log$-concave on $\R$. It follows that the right derivative of $F$, denoted by $F'_+$, exists everywhere and that $t \mapsto t F'_+(t) / F(t) $ is non-increasing on $(0, + \infty)$.

Notice that the function $F$ is $\frac{1}{n}$-concave on $\R_+$ if and only if the function
$$ t \mapsto \frac{t F'_+(t)}{F(t)} \frac{F(t)^{\frac{1}{n}}}{t} $$
is non-increasing on $\R_+$. A direct change of variables shows that
$$ \frac{F(t)}{t^n} = \int_A \phi(tx) \, \de x. $$
By assumption, the function $t \mapsto \phi(tx) $ is non-increasing on $\R_+$. It follows that the function $t \mapsto F(t)^{\frac{1}{n}} / t$ is non-increasing on $(0, + \infty)$. Hence, the function $t \mapsto (F(t)^{\frac{1}{n}})'_+$ is non-increasing on $(0, + \infty)$ as the product of two non-negative non-increasing functions on $(0, + \infty)$. We conclude that $F$ is $\frac{1}{n}$-concave on $\R_+$.
\end{proof}

\paragraph*{Remarks.}

\begin{enumerate}
\item Proposition~\ref{equiv} is established in~\cite{GZ} for the Gaussian measure by differentiating twice.
\item The assumption {\it $t \mapsto \phi(tx) $ is non-increasing on $\R_+$} is satisfied if $\phi$ is an even $-\infty$-concave function.
\item The converse of Proposition~\ref{equiv} is false in general, by taking, for example, $\phi(x) = 1_{[-1, 2]}(x)$, $x \in \R$, and $A = [-1,1]$.
\end{enumerate}

Proposition~\ref{equiv} is related to the (B)-conjecture. This conjecture was posed by W. Banaszczyk~\cite{Latala} and asks whether the function $t \mapsto \gamma_n(\e^t A)$ is $\log$-concave on $\R$, for every convex symmetric set $A \subset \R^n$. The (B)-conjecture was proved by Cordero-Erausquin, Fradelizi and Maurey in~\cite{CFM}. In the same paper~\cite{CFM}, the authors have also shown that for every unconditional $\log$-concave measure $\mu$ on $\R^n$ and for every unconditional convex subset $A \subset \R^n$, the function $t \mapsto \mu(\e^t A)$ is $\log$-concave on $\R$. Using this and the point 2. of the previous remark, we may apply Proposition~\ref{equiv} to obtain the following corollary:

\begin{coro}\label{dilate}

Let $\mu$ be an unconditional $\log$-concave measure on $\R^n$ and let $A$ be an unconditional convex subset of $\R^n$. Then, the measure $\mu$ is $\frac{1}{n}$-concave in the class of dilates of $A$. More precisely, for every $A_1, A_2 \in \{\alpha A; \alpha > 0 \}$ and for every $\lambda \in [0,1]$, we obtain
$$ \mu((1-\lambda)A_1 + \lambda A_2)^{\frac{1}{n}} \geq (1-\lambda) \mu(A_1)^{\frac{1}{n}} + \lambda \mu(A_2)^{\frac{1}{n}}. $$

\end{coro}

\paragraph*{Remark.}

Very recently, Livne Bar-on~\cite{L} and Saroglou~\cite{S} proved, using different methods, that in dimension~2 for the uniform measure $\mu_K$ on a symmetric convex set $K \subset \R^2$ ($i.e.$ $\de \mu_K(x) = 1_K(x) \, \de x$), the function $t \mapsto \mu_K(\e^t A)$ is log-concave on $\R$ for every symmetric convex set $A \subset \R^2$. However, for our problem, this information is not useful since the uniform measure on a convex subset of $\R^n$ is a $\frac{1}{n}$-concave measure. \\

A natural question is to ask if the role of the symmetry can be relaxed. When restricted to the Gaussian measure, it has been shown by Nayar and Tkocz in \cite{N-Tkocz}, that in dimension~2 there exist non-symmetric convex sets $A, B \subset \R^2$ satisfying $0 \in A \subset B$ and
\begin{eqnarray}\label{CE}
\gamma_2 \left(\frac{A+B}{2}\right)^{\frac{1}{2}} < \frac{ \gamma_2(A)^{\frac{1}{2}} + \gamma_2(B)^{\frac{1}{2}}}{2}.
\end{eqnarray}
One can then construct an explicit counterexample in every dimension $n \geq 2$. Moreover, the counterexample in~\cite{N-Tkocz} shows more than inequality~(\ref{CE}). It shows that
\begin{eqnarray}\label{CE2}
\gamma_2 \left(\frac{A+B}{2}\right)^s < \frac{ \gamma_2(A)^s + \gamma_2(B)^s}{2},
\end{eqnarray}
for every $s \geq 1 - \frac{2}{\pi}$. However, it could be of interest to know what happens when $s \in (0, 1 - \frac{2}{\pi})$.

Notice that the same counterexample with the following log-concave unconditional measure instead of the Gaussian measure
$$ \de \mu(x,y) = \e^{-|x|}\e^{-|y|} \, \de x \, \de y, \quad (x,y) \in \R^2 $$
satisfies inequality~(\ref{CE2}) for every $s>0$.

Thus, in general, the symmetry assumption on the measure is not sufficient. \\

On the other hand, the concavity of a non-symmetric convex measure cannot be improved in general in the class of symmetric sets even in dimension~1:

\begin{prop}

Let $0 < s < 1$ and $r > s$. There exists an $s$-concave measure $\mu$ on $\R$ and symmetric sets $A, B \subset \R$ such that
$$ \mu\left(\frac{A+B}{2}\right) < \left( \frac{\mu(A)^r + \mu(B)^r}{2} \right)^{\frac{1}{r}}. $$

\end{prop}

\begin{proof}
Let us define $ \de \mu (x) = x^{1/\gamma}1_{\R_+}(x) \, \de x$, with $\gamma=\frac{s}{1-s} > 0$.
Let us consider the sets $A=[-a,a]$ and $B=[-b,b]$ with $0<a<b$. Notice that
$$ \lim_{a \to 0} \mu\left(\frac{A+B}{2}\right) = {\mu \left(\frac{B}{2} \right)} = \frac{\mu(B)}{2^{\frac{1}{s}}} = \lim_{a \to 0} \left( \frac{\mu(A)^s + \mu(B)^s}{2} \right)^{\frac{1}{s}}. $$
Since $\mu(A) \neq \mu(B)$, it follows by comparing means that
$$ \left( \frac{\mu(A)^s + \mu(B)^s}{2} \right)^{\frac{1}{s}} < \left( \frac{\mu(A)^r + \mu(B)^r}{2} \right)^{\frac{1}{r}}. $$
we conclude that for sufficiently small $a$,
$$ \mu\left(\frac{A+B}{2}\right) < \left( \frac{\mu(A)^r + \mu(B)^r}{2} \right)^{\frac{1}{r}}. $$
\end{proof}

Thus, in general, the symmetry assumption on the sets is not sufficient.

\subsection{The case of sets with a maximal section of equal measure}

In this section, we consider $\mathcal{C}$ to be the class of couples of Borel subsets of $\R^n$ having a maximal section of equal measure. A famous result of Bonnesen~\cite{BF} (for convex sets) states that if $A,B \subset \R^n$ satisfy
$$ \sup_{t \in \R} |A \cap (u^{\perp} + tu)|_{n-1} = \sup_{t \in \R} |B \cap (u^{\perp} + tu)|_{n-1}, $$
for a certain $u \in \mathcal{S}^{n-1}$, where $|\cdot|_{n-1}$ denotes the $(n-1)$-dimensional Lebesgue measure, then for every $\lambda \in [0,1]$, one has
$$ |(1-\lambda) A + \lambda B| \geq (1-\lambda) |A| + \lambda |B|. $$
There exists a functional version of Bonnesen's result established by Henstock and Macbeath~\cite{HMc} in dimension~1 (see Proposition~\ref{function} below) and later on by Dancs and Uhrin~\cite{DU} in higher dimension (see Proposition~\ref{function-n} below).

\begin{prop}[Henstock, Macbeath~\cite{HMc}]\label{function}

Let $\lambda \in [0,1]$. Let $f,g,h : \R \to \R_+$ be non-negative measurable functions such that $\max(f) = \max(g)$ and such that for every $x,y \in \R$
$$ h((1-\lambda)x + \lambda y) \geq \min(f(x),g(y)). $$
Then, one has
$$ \int_{\R} h(x) \, \de x \geq (1-\lambda) \int_{\R} f(x) \, \de x + \lambda \int_{\R} g(x) \, \de x. $$

\end{prop}

We deduce the following result:

\begin{prop}\label{concave}

Let $\phi : \R \to \R_+$ be a $-\infty$-concave function such that $\max(\phi)$ is attained at $a \in \R$. Let $\mu$ be a measure with density function $\phi$. Let $A,B$ be Borel subsets of $\R$ such that $a \in A \cap B$. Then, for every $\lambda \in [0,1]$, we have
$$ \mu((1-\lambda) A + \lambda B) \geq (1-\lambda) \mu(A) + \lambda \mu(B). $$

\end{prop}

\begin{proof}
Let $\lambda \in [0,1]$. We define, for every $x \in \R$, $h(x)=\phi(x) 1_{(1-\lambda)A + \lambda B}(x)$, $f(x)=\phi(x) 1_A(x)$, $g(x)=\phi(x) 1_B(x)$. Notice that for every $x,y \in \R$ one has
$$ h((1-\lambda)x + \lambda y) \geq \min(f(x),g(y)), $$
and $\max(f) = \max(g) = \phi(a)$. It follows from Proposition~\ref{function} that
$$ \int_{\R} h(x) \, \de x \geq (1-\lambda) \int_{\R} f(x) \, \de x + \lambda \int_{\R} g(x) \, \de x. $$
In other words, we obtain
$$ \mu((1-\lambda) A + \lambda B) \geq (1-\lambda) \mu(A) + \lambda \mu(B). $$
\end{proof}

\paragraph*{Remark.}

Proposition~\ref{concave} was established in~\cite{GZ} for the case where $\mu$ is the Gaussian measure on $\R$ and where $A, B \subset \R$ are convex. In the same paper, the authors were able to remove the convexity assumption for only one set, by using long computations. Our method bypasses the use of geometric tools and relies on the functional version Proposition~\ref{function}. \\

Conversely, if a measure $\mu$ on $\R$, with density function $\phi$ with respect to the Lebesgue measure, satisfies
$$ \mu((1-\lambda)A + \lambda B) \geq (1-\lambda)\mu(A) + \lambda \mu(B), $$
for every $\lambda \in [0,1]$ and for all symmetric convex sets $A,B \subset \R$, then one has for every $\lambda \in [0,1]$ and for every $a,b \in \R_+$,
$$ \int_{-((1-\lambda)a + \lambda b)}^{(1-\lambda)a + \lambda b} \phi(x) \, \de x \geq (1-\lambda) \int_{-a}^a \phi(x) \, \de x + \lambda \int_{-b}^b \phi(x) \, \de x. $$
It follows that the function $t \mapsto \phi(t) + \phi(-t)$ is non-increasing on $\R_+$. Notice that this condition is satisfied for more general functions than $-\infty$-concave functions attaining the maximum at $0$.

However, one can use the same argument to see that if one assumes $A,B \subset \R$ convex containing~$0$ (not necessarily symmetric), then it follows that the density function $\phi$ is necessarily non-decreasing on $(-\infty; 0]$ and non-increasing on $[0;+\infty)$. Notice that this is equivalent to the fact that the density function $\phi$ is $-\infty$-concave and $\max(\phi)$ is attained at~$0$. \\

By tensorization, Proposition~\ref{concave} leads to the following corollary:

\begin{coro}\label{n-concave}

Let $\mu_i$, $1 \leq i \leq n$, be measures with densities $\phi_i : \R \to \R_+$ such that $\phi_i$ are non-decreasing on $(-\infty; 0]$ and non-increasing on $[0; + \infty)$. Let $\mu$ be the product measure of $\mu_1, \dots, \mu_n$ and let $A,B \subset \R^n$ be the product of $n$ Borel subsets of $\R$ such that $0 \in A \cap B$. Then, for every $\lambda \in [0,1]$, we have
$$ \mu((1-\lambda) A + \lambda B)^{\frac{1}{n}} \geq (1-\lambda) \mu(A)^{\frac{1}{n}} + \lambda \mu(B)^{\frac{1}{n}}. $$
\end{coro}

\begin{proof}
We follow~\cite{GZ}. By assumption, $A = \Pi_{i=1}^n A_i$ and $B = \Pi_{i=1}^n B_i$, where for every $i \in \{1, \dots, n\}$, $A_i$ and $B_i$ are Borel subsets of $\R$ containing $0$. Let $\lambda \in [0,1]$. Notice that
$$ (1-\lambda)A + \lambda B = \Pi_{i=1}^n ((1-\lambda)A_i + \lambda B_i). $$
Using Proposition~\ref{concave} and an inequality of Minkowski (see $e.g.$ \cite{HLP}), one deduces that
\begin{eqnarray*}
\mu((1-\lambda) A + \lambda B)^{\frac{1}{n}} & = & \left( \Pi_{i=1}^n \,  \mu_i((1-\lambda)A_i + \lambda B_i) \right)^{\frac{1}{n}} \\ & \geq & \left( \Pi_{i=1}^n \, ((1-\lambda) \mu_i(A_i) + \lambda \mu_i(B_i)) \right)^{\frac{1}{n}} \\ & \geq & \left( \Pi_{i=1}^n (1-\lambda) \mu_i(A_i) \right)^{\frac{1}{n}} + \left( \Pi_{i=1}^n \lambda \mu_i(B_i) \right)^{\frac{1}{n}} \\ & = & (1-\lambda) \mu(A)^{\frac{1}{n}} + \lambda \mu(B)^{\frac{1}{n}}.
\end{eqnarray*}
\end{proof}

Another consequence of Proposition~\ref{concave} is that certain particular product measures are concave measures if $A$ is a union of parallel slabs containing the origin.

\begin{coro}\label{slab}

Let $\mu_1$ be a measure with density function $\phi: \R \to \R_+$, such that $\phi$ is non-decreasing on $(-\infty; 0]$ and non-increasing on $[0; + \infty)$. Let $\mu_2$ be a $(n-1)$-dimensional measure and let $\mu$ be the product measure of $\mu_1$ and $\mu_2$. Let $A = A_1 \times \R^{n-1}$, where $A_1$ is a Borel subset of $\R$ and let $B$ be a Borel subset of $\R^n$ such that $0 \in A \cap B$. Then, for every $\lambda \in [0,1]$, we have
\begin{eqnarray}\label{unionslab}
\mu((1-\lambda) A + \lambda B) \geq (1-\lambda) \mu(A) + \lambda \mu(B).
\end{eqnarray}

\end{coro}

Inequality~(\ref{unionslab}) was established in~\cite{GZ} with power $1/n$ for slabs, for the case where $\mu$ is the Gaussian measure.

\begin{proof}
We follow~\cite{GZ}. Let us denote $B_S = \proj_{e_1}(B) \times \R^{n-1}$, where $\proj_{e_1}(B)$ denotes the orthogonal projection of $B$ onto the first coordinate axis. Then, for every $\lambda \in [0,1)$, one has
$$ (1-\lambda) A + \lambda B = (1-\lambda) A + \lambda B_S. $$
It follows, using Proposition~\ref{concave}, that
\begin{eqnarray*}
\mu((1-\lambda) A + \lambda B) & = & \mu((1-\lambda) A + \lambda B_S) \\ & = & \mu(((1-\lambda) A_1 + \lambda \proj_{e_1}(B)) \times \R^{n-1}) \\ & = & \mu_1((1-\lambda) A_1 + \lambda \proj_{e_1}(B)) \cdot \mu_2(\R^{n-1}) \\ & \geq & ((1-\lambda) \mu_1(A_1) + \lambda \mu_1(\proj_{e_1}(B))) \cdot \mu_2(\R^{n-1}) \\ & = & (1-\lambda) \mu(A) + \lambda \mu(B_S) \\ & \geq & (1-\lambda) \mu(A) + \lambda \mu(B).
\end{eqnarray*}
\end{proof}

On the other hand, Proposition~\ref{function} can be generalized to the dimension~$n$. First, let us define for a non-negative measurable function $f : \R^n \to \R_+$ and for $u \in \mathcal{S}^{n-1}$,
$$ m_u(f) = \sup_{t \in \R} \int_{u^{\perp}} f(x+tu) \, \de x. $$

\begin{prop}[Dancs, Uhrin \cite{DU}]\label{function-n}

Let $-\frac{1}{n-1} \leq \gamma \leq + \infty$, $\lambda \in [0,1]$ and $f,g,h : \R^n \to \R_+$ be non-negative measurable functions such that for every $x,y \in \R^n$,
$$ h((1-\lambda)x + \lambda y) \geq \left( (1-\lambda)f(x)^{\gamma} + \lambda g(y)^{\gamma} \right)^{\frac{1}{\gamma}}. $$
If there exists $u \in \mathcal{S}^{n-1}$ such that $m_u(f) = m_u(g) < +\infty$, then
$$ \int_{\R^n} h(x) \, \de x \geq (1-\lambda) \int_{\R^n} f(x) \, \de x + \lambda \int_{\R^n} g(x) \, \de x. $$

\end{prop}

Let us denote for a measure $\mu$ with density function $\phi$, for a Borel subset $A \subset \R^n$ and for a hyperplane $H \subset \R^n$,
$$ \mu_{n-1}(A \cap H) = \int_{A \cap H} \phi(x) \, \de x. $$
We deduce the following result:

\begin{prop}\label{rouge}

Let $\mu$ be a measure with density function $\phi : \R^n \to \R_+$ such that $\phi$ is $-\frac{1}{n-1}$-concave. Let $A,B$ be Borel subsets of $\R^n$. If there exists $u \in \mathcal{S}^{n-1}$ such that
$$ \sup_{t \in \R} \mu_{n-1}(A \cap (u^{\perp} + tu)) = \sup_{t \in \R} \mu_{n-1}(B \cap (u^{\perp} + tu)), $$
then, for every $\lambda \in [0,1]$, we have
$$ \mu((1-\lambda) A + \lambda B) \geq (1-\lambda) \mu(A) + \lambda \mu(B). $$

\end{prop}

\begin{proof}
Let $\lambda \in [0,1]$. Let us take $f = \phi 1_A$, $g=\phi 1_B$ and $h=\phi 1_{(1-\lambda)A + \lambda B}$. Then, for every $x, y \in \R^n$, one has
$$ h((1-\lambda)x + \lambda y) \geq \left( (1-\lambda)f(x)^{\gamma} + \lambda g(y)^{\gamma} \right)^{\frac{1}{\gamma}}, $$
where $\gamma = -\frac{1}{n-1}$. Moreover,
$$ \int_{u^{\perp}} f(x+tu) \, \de x = \int_{A \cap (u^{\perp} + tu)} \phi(x) \, \de x = \mu_{n-1}(A \cap (u^{\perp} + tu)). $$
It follows that $m_u(f) = m_u(g)$. From Proposition~\ref{function-n}, we obtain that
$$ \mu((1-\lambda) A + \lambda B) \geq (1-\lambda) \mu(A) + \lambda \mu(B). $$
\end{proof}

\section{Application to the parallel volume}

Let us see how improvements of the concavity of sub-convex measures can improve the concavity of a generalized form of the parallel volume. The parallel volume of a Borel subset $A$ of $\R^n$, namely the function $t \mapsto |A + t B_2^n|$, is a particularly useful function in Geometry, which has been highlighted by the precursor works of Steiner~\cite{Steiner}. For more modern applications, the parallel volume and its generalized forms are still studied (see $e.g.$ \cite{HLW}, \cite{Kampf}). Moreover, this notion of parallel volume leads to the powerful theory of mixed volumes (see~\cite{Schneider} for further details). 

As a consequence of the Brunn-Minkowski inequality~(\ref{BM}), one can see that if $A \subset \R^n$ is convex, then the parallel volume of $A$ is $\frac{1}{n}$-concave on $\R_+$. More generally, if a measure $\mu$ is $s$-concave, with $s \in [-\infty;+\infty]$, in the class of sets of the form $\{ A+tB; t \in \R_+\}$, where $A$ and $B$ are convex subsets of $\R^n$, then the generalized parallel volume $t \mapsto \mu(A + tB)$ is $s$-concave on $\R_+$. Indeed, for every $t_1, t_2 \in \R_+$ and for every $\lambda \in [0,1]$, one has
\begin{eqnarray*}
\mu(A + ((1-\lambda)t_1 + \lambda t_2)B) & = & \mu((1-\lambda)(A+t_1B) + \lambda (A+t_2 B)) \\ & \geq & \left( (1-\lambda) \mu(A + t_1 B)^s + \lambda \mu(A + t_2 B)^s \right)^{\frac{1}{s}}.
\end{eqnarray*}
Using this and Corollary~\ref{n-concave}, we obtain the following corollary:

\begin{coro}

Let $\mu_i$, $1 \leq i \leq n$, be measures with densities $\phi_i : \R \to \R_+$ such that $\phi_i$ are non-decreasing on $(-\infty; 0]$ and non-increasing on $[0; + \infty)$. Let $\mu$ be the product measure of $\mu_1, \dots, \mu_n$ and let $A,B \subset \R^n$ be coordinate boxes containing the origin. Then the function $ t \mapsto \mu(A + t B) $ is $\frac{1}{n}$-concave on $\R_+$.

\end{coro}

In the case of non-convex sets, this concavity property is false in general, even for the classical parallel volume $|A+t B_2^n|$. However, some conditions are given on $A$ in~\cite{FM} for which the parallel volume of $A$ is $\frac{1}{n}$-concave on $\R_+$. Notice that other concavity properties of generalized forms of the classical parallel volume have been established in~\cite{M}.

\bibliographystyle{amsplain}

\begin{thebibliography}{30}


\bibitem{BF} T. Bonnesen, W. Fenchel, {\it Theorie der konvexen K\"orper}, Springer, Berlin, 1934. English translation: {\it Theory of convex bodies}, edited by L. Boron, C. Christenson and B. Smith. BCS Associates, Moscow, ID, 1987.
	\bibitem{Borell} C. Borell, {\it Convex set functions in d-space}, Periodica Mathematica Hungarica Vol. 6, 111-136, 1975.
	\bibitem{BL} H. J. Brascamp, E. H. Lieb, {\it On extensions of the Brunn-Minkowski and Pr\'ekopa-Leindler theorems, including inequalities for log concave functions, and with an application to the diffusion equation}, (1976) J. Funct. Anal. 22 366-389.
	\bibitem{CFM} D. Cordero-Erausquin, M. Fradelizi, B. Maurey, {\it The (B) conjecture for the Gaussian measure of dilates of symmetric convex sets and related problems}, J. Funct. Anal. 214 (2004), no. 2, 410-427.
	\bibitem{DU} S. Dancs, B. Uhrin, {\it On a class of integral inequalities and their measure-theoretic consequences}, J. Math. Anal. Appl. 74 (1980), no. 2, 388-400.
	\bibitem{FM} M. Fradelizi, A. Marsiglietti, {\it On the analogue of the concavity of entropy power in the Brunn-Minkowski theory}, Advances in Applied Mathematics 57 (2014), 1-20.
	\bibitem{Gardner} R. J. Gardner, {\it The Brunn-Minkowski inequality}, Bull. Amer. Math. Soc. (N.S.) 39 (2002), no. 3, 355-405.
	\bibitem{GZ} R. J. Gardner, A. Zvavitch, {\it Gaussian Brunn-Minkowski inequalities}, Trans. Amer. Math. Soc. 362 (2010), no. 10, 5333-5353.
	\bibitem{HLP} G. H. Hardy, J. E. Littlewood, G. P\'olya, {\it Inequalities}, Cambridge University Press, Cambridge, 1959.
	\bibitem{HMc} R. Henstock, A. M. Macbeath, {\it On the measure of sum-sets. I. The theorems of Brunn, Minkowski, and Lusternik}, Proc. London Math. Soc. (3) 3, (1953). 182-194.
	\bibitem{HN} M. A. Hern\'andez Cifre, J. Yepes Nicol\'as, {\it Refinements of the Brunn-Minkowski inequality}, J. Convex Anal. 21 (3) (2014), 1-17.
	\bibitem{HLW} D. Hug, G. Last, W. Weil, {\it A local Steiner-type formula for general closed sets and applications}, Math. Z. 246 (2004), no. 1-2, 237-272.
	\bibitem{Kampf} J. Kampf, {\it The parallel volume at large distances}, Geom. Dedicata 160 (2012), 47-70.
	\bibitem{Latala} R. Lata\l a, {\it On some inequalities for Gaussian measures}, (English summary) Proceedings of the International Congress of Mathematicians, Vol. II (Beijing, 2002), 813-822, Higher Ed. Press, Beijing, 2002.
	\bibitem{Leindler} L. Leindler, {\it On a certain converse of H\"older's inequality}, II, Acta Sci. Math., 33 (1972), 217-223.
	\bibitem{L} A. Livne Bar-on, {\it The (B) conjecture for uniform measures in the plane}, preprint, arXiv:1311.6584 [math.FA].
	\bibitem{M} A. Marsiglietti, {\it Concavity properties of extensions of the parallel volume}, To appear in Mathematika.
	\bibitem{N-Tkocz} P. Nayar, T. Tkocz, {\it A note on a Brunn-Minkowski inequality for the Gaussian measure}, Proc. Amer. Math. Soc. 141 (2013), no. 11, 4027-4030.
	\bibitem{Prekopa71} A. Pr\'ekopa, {\it Logarithmic concave measures with application to stochastic programming}, Acta Sci. Math., 32 (1971), 301-316.
	\bibitem{Prekopa} A. Pr\'ekopa, {\it On logarithmic concave measures and functions}, Acta Sci. Math. (Szeged) 34 (1973), 335-343.
	\bibitem{S} C. Saroglou, {\it Remarks on the conjectured log-Brunn-Minkowski inequality}, To appear in Geom. Dedicata.
	\bibitem{Schneider} R. Schneider, {\it Convex bodies: the Brunn-Minkowski theory}, Encyclopedia of Mathematics and its Applications, 44. Cambridge University Press, Cambridge, 1993. xiv+490 pp. ISBN: 0-521-35220-7.
	\bibitem{Steiner} J. Steiner, {\it \"Uber parallele Fl\"achen}, Monatsbericht der Akademie der Wissenschaften zu Berlin (1840), pp. 114-118.
\end{thebibliography}

\end{document}